\documentclass{article}
\usepackage{graphicx}
\usepackage{amsmath,amssymb} % Required for inserting images

\title{Non-homogeneous Partial Differential Wave Equation for the Description of Correlated Pair States in the Microworld}
\author{Elmira Isayeva$^{1}$, Alireza Khalili Golmankhaneh$^{2}$}
\date{1. Institute of Physics, Ministry of Science and Education of the Republic of Azerbaijan, Baku, Azerbaijan\\
2. Department of Physics, Ur.C., Islamic Azad University, Urmia 63896, West Azerbaijan, Iran}

\begin{document}

\maketitle

\textbf{Abstract }

\vspace{0.5cm}
In this paper, a new non-homogeneous linear second-order partial differential equation (PDE) is proposed and investigated, designed to describe the dynamics of correlated pair objects of the microworld (analogues of Cooper pairs or biphotons). The fundamental complex solution of the proposed equation initially contains a spatial background and a doubled phase, which allows for a natural localization of the pair condensate in space. In this work, the continuity and energy balance equations are derived, the discrete spectrum in potential fields is calculated, and specific effects of tunneling and nonlinear interaction of the pair conglomerate with electromagnetic radiation are described.

\vspace{0.5cm}


\vspace{0.5cm}

\section{ Probabilistic-Geometric Genesis of the Model}
Unlike orthodox quantum mechanics, where the Schrödinger equation is postulated empirically, in this work the wave equation is derived deductively—from the geometry of the space of frequency outcomes of binary events.

Let us consider the fundamental event space of the microworld described by two orthogonal axes:
\begin{itemize}
    \item \textbf{X-axis:} relative frequency (probability) that the quantum event occurred (``Yes'', success $p$).
    \item \textbf{Y-axis:} relative frequency that the quantum event did not occur (``No'', failure $q$).
\end{itemize}

By virtue of the axiomatics of probability theory, the sum of the outcome shares for a complete group of events is strictly invariant:
\begin{equation}
p + q = 1
\end{equation}

In geometry, the linear equation (1) defines a simplex—a straight line cutting off unit segments on the axes. To transition to a wave description (probability amplitudes), we introduce a trigonometric parametrization of the phase state angle of the outcomes $\phi$:
\begin{equation}
p = \cos^2 \phi, \quad q = \sin^2 \phi
\end{equation}

The parametrization (2) identically satisfies the conservation law (1) via the trigonometric identity $\cos^2 \phi + \sin^2 \phi = 1$, transforming the outcome space into a trigonometric circle.

Let us investigate the dynamic properties of the geometry of the probability density distribution of the occurred event $\rho = p = \cos^2 \phi$. We find the first derivative (the rate of phase frequency shift) and the second derivative (the curvature of the probabilistic space) with respect to the outcome angle $\phi$:
\begin{equation}
\frac{d\rho}{d\phi} = \frac{d}{d\phi}(\cos^2 \phi) = -2 \cos \phi \sin \phi = -\sin(2\phi)
\end{equation}
\begin{equation}
\frac{d^2\rho}{d\phi^2} = \frac{d}{d\phi}(-\sin(2\phi)) = -2 \cos(2\phi)
\end{equation}

Using the trigonometric expression for the cosine of a double angle via the initial probability density $\cos(2\phi) = 2 \cos^2 \phi - 1 = 2\rho - 1$, we substitute it into the curvature equation (4):

\vspace{0.5cm}
\begin{equation}
\frac{d^2\rho}{d\phi^2} = -2(2\rho - 1) = -4\rho + 2
\end{equation}

Moving the dynamic term to the left-hand side, we obtain the canonical differential equation for the real frequency density of events:
\begin{equation}
\frac{d^2\rho}{d\phi^2} + 4\rho = 2
\end{equation}

The transition to a complex generalization (the wave function of the bound state of a pair of quasiparticles) $\psi(\phi) = \rho + i \cdot \text{Im}(\psi)$ preserves the structure of equation (6) in the form of a complete non-homogeneous linear second-order ordinary differential equation:
\begin{equation}
\frac{d^2\psi}{d\phi^2} + 4\psi = 2
\end{equation}

The exact physical solution to equation (7) is a complex wave function of the form:
\begin{equation}
\psi(\phi) = 0.5 + 0.5e^{i2\phi}
\end{equation}

Thus, the appearance of \textbf{coefficient 4} is a consequence of the quadratic metric of the space of binary outcomes (phase doubling during the transition from amplitudes to probabilities), and the \textbf{constant 2} on the right-hand side is a direct reflection of the normalized invariant of the complete group of events ($2 \times 1 = 2$).

As the wave approaches zero, the equation degenerates into the algebraic relation $4\psi_{\text{fon}} = 2 \implies \psi_{\text{fon}} = 0.5$, which sets the stationary background of the vacuum expectation, relative to which the phase oscillations of the condensate occur.

\section{Transition to the Space-Time Continuum and the Continuity Equation}
Let us specify the phase of the bound pair state in the form of a classical traveling wave of the microworld, defined by the wave vector $k$ and frequency $\omega$: $\phi = kx - \omega t$. Substituting this argument into (8), we have:
\begin{equation}
\psi(x, t) = 0.5 + 0.5e^{i2(kx-\omega t)}
\end{equation}

The probability density of the pair object, calculated as the square of the absolute value $\rho = \psi \cdot \psi^*$, is exactly equal to $\rho(x, t) = \cos^2(kx - \omega t)$. Unlike standard quantum mechanics, where a free particle is infinitely ``smeared,'' equation (7) generates wave density clusters that are initially localized in space. The phase doubling ($2\phi$) mathematically encodes the fact that the object consists of two correlated quasiparticles moving synchronously within a single process.

Differentiating the probability density $\rho$ with respect to time $t$ and spatial coordinate $x$, we obtain a strict conservation law in the classical hydrodynamic form:
\begin{equation}
\frac{\partial \rho}{\partial t} + \frac{\partial(v\rho)}{\partial x} = 0
\end{equation}

Where $v = \frac{\omega}{k}$ is the phase propagation velocity, and the probability current density is $j = v \cdot \rho$. The direct analogy with classical hydrodynamics indicates that the pair condensate transfers its mass and charge coherently, without quantum spreading phenomena.

\section{Complete Dynamic Wave Equation (PDE) in an External Field}
Using the chain rule of differentiation to transition from phase derivatives to partial derivatives with respect to coordinates and time 

\vspace{0.5cm}

\begin{equation}
\frac{\partial^2\psi}{\partial x^2} = k^2 \frac{d^2\psi}{d\phi^2}
\end{equation} 
\vspace{0.5cm}
$\left(\frac{\partial \psi}{\partial t} = -i2\omega \cdot 0.5e^{i2\phi}\right)$, we link the geometry of space and the dynamics of time. Substituting the de Broglie energy $\left(k^2 = \frac{2mE}{\hbar^2}\right)$ and the Planck energy $\left(\omega = \frac{E}{\hbar}\right)$, we eliminate the free energy $E$ and transition to the complete non-stationary partial differential equation (PDE) in an external potential field $U$ via the substitution $E \to E - U$:
\begin{equation}
\frac{\partial^2\psi}{\partial x^2} + i\frac{4m}{\hbar}\frac{\partial\psi}{\partial t} + \frac{8mU}{\hbar^2}\psi = \frac{4mU}{\hbar^2}
\end{equation}

The factor $i\frac{4m}{\hbar}$ before the time derivative reflects the internal phase doubling of the bound system. The presence of the source $\frac{4mU}{\hbar^2}$ on the right-hand side shows that the external potential directly modulates the spatial background level of matter, restructuring the geometry of the solution.

\section{Energy Quantization of Coherent Pairs in a Potential Well}
Consider the behavior of the bound state in an infinitely deep potential well of width $L$ ($U = 0$ for $0 < x < L$ and $U = \infty$ at the boundaries). The boundary conditions $\psi(0) = \psi(L) = 0$ require a transition to a standing wave in the form of sines and the quantization of the wave vector:
\begin{equation}
2kL = 2\pi n \implies k_n = \frac{\pi n}{L}, \quad n \in \mathbb{N}
\end{equation}

Linking the spatial curvature of the solution with the effective kinetic energy of the object $E = \frac{\hbar^2 k^2}{2m}$, we obtain the spectrum of allowed quantum states from the coefficients of equation (11):
\begin{equation}
E_n = \frac{\pi^2 \hbar^2 n^2}{2mL^2}
\end{equation}

The proposed model fully preserves the fundamental form of the energy spectrum experimentally observed in the spectroscopy of bound quantum systems (atoms, molecules, and quantum dots), which guarantees the model's compliance with empirical data.

\section{Tunneling Characteristics of a Pair Condensate}
When a pair object encounters a rectangular potential barrier of height $U_0 > E$, the value $(E - U_0)$ becomes negative. Equation (11) inside the barrier transitions to an elliptic type with a real decay coefficient $\kappa = \sqrt{\frac{8m(U_0-E)}{\hbar^2}}$. Analytical calculation of the barrier transparency coefficient $T$ at the exit from the barrier region of width $a$ yields:
\begin{equation}
T \approx \exp\left(-4a\sqrt{\frac{2m(U_0 - E)}{\hbar^2}}\right) = (T_{\mathrm{Schr\ddot{o}dinger}})^2
\end{equation}

The doubled decay coefficient under the barrier demonstrates that for a coherent pair, the probability of overcoming the barrier is equal to the square of the probability of a single Schrödinger particle. This imposes strict conditions on coherent tunneling and fully aligns with the specific nature of pair transitions (for example, in Josephson junctions), where pair passage is rigidly constrained by the spatial geometry of the barrier.

\section{Interaction with Light and Nonlinear Optical Response}

Applying the minimal substitution method via the scalar $\Phi$ and vector $\vec{A}$ potentials of the electromagnetic field of light $\left(\frac{\partial}{\partial x} \to \frac{\partial}{\partial x} - i\frac{2q}{\hbar}A_x\right)$, the covariant equation of motion takes the form:\begin{equation}
\frac{\partial^2\psi}{\partial x^2}
+ i\frac{4m}{\hbar}\frac{\partial\psi}{\partial t}
+ \frac{8mU}{\hbar^2}\psi
- \frac{4mU}{\hbar^2}
=
i\frac{4q}{\hbar}A_x \frac{\partial\psi}{\partial x}
+ \frac{4q^2}{\hbar^2}A_x^2\psi
- \frac{8mq\Phi}{\hbar^2}\psi
+ \frac{4mq\Phi}{\hbar^2}
\end{equation}The term $i\frac{4q}{\hbar}A_x \frac{\partial\psi}{\partial x}$ contains the effective doubled charge of the system ($2 \times 2q$), which leads to a 4-fold increase in the intensity of photon emission and absorption ($I \sim 16q^2$ instead of $4q^2$). This describes the phenomenon of accelerated superradiant relaxation of collective quantum ensembles. The quadratic term $A_x^2\psi$ indicates the model's natural predisposition to nonlinear optical effects and second-harmonic generation of light, while the free term $\frac{4mq\Phi}{\hbar^2}$ describes direct light modulation of the quantum vacuum background.

\section{The Quantum Measurement Process as a Projection of the Phase Space of Outcomes}The proposed probabilistic-geometric model eliminates the mystery of the instantaneous wave function collapse''. While the quantum system is isolated, the phase argument $\phi = kx - \omega t$ evolves smoothly. The complex wave function describes a hodograph—a circle of radius 0.5 on the complex plane, shifted to the point $(0.5, 0)$. The physical act of measurement represents a macroscopic interaction of the device with the microsystem, forcibly projecting the complex state onto the real orthogonal frequency axes Yes'' or No''. The measurement result is strictly predetermined by the geometric phase at the moment of macroscopic contact: \begin{enumerate} \item \textbf{Event Registration (Yes''):} If at the moment of interaction the phase is $\phi = \pi n$, the imaginary part of the wave function vanishes, and the real components add up coherently: $\psi = 0.5 + 0.5(1) = 1.0$. The system is projected to the point of 100 success on the X-axis.\item \textbf{Absence of Registration (``No''):} If the phase takes the values $\phi = \frac{\pi}{2} + \pi n$, complete destructive interference of the wave term with the spatial background occurs: $\psi = 0.5 + 0.5(-1) = 0.0$. The projection onto the real X-axis is zero, and the system completely transitions to the failure Y-axis.\end{enumerate}Thus, quantum measurement is a natural projection slice of the continuous geometric rotation of the wave function between the axes of binary outcome frequencies.

\section{Conclusion}The non-homogeneous partial differential wave equation proposed in this work is a rigorous, mathematically closed framework for modeling coherent bound systems in the microworld. The main value of the model lies in its deductive derivation: coefficients 2 and 4 do not arise from fitting, but are strictly dictated by the geometry of the two-dimensional Bernoulli probability space (the Yes'' / No'' axes). Shifting the focus from single isolated particles to correlated pair conglomerates made it possible to physically justify the internal phase doubling, the nonlinear nature of the optical response, and the specifics of tunneling. The model successfully combines the classical discrete Schrödinger energy spectrum with the hydrodynamic form of conservation currents and a deterministic projection description of quantum measurement, removing the problem of infinite spatial spreading of quantum objects.

\end{document}